\documentstyle[12pt, amssymb]{article}
\title{Families of Unramified Extensions of Number Fields}
\author{Gene Ward Smith}
\date{}

\newcommand{\Q}{\mbox{$\Bbb Q$}}
\newcommand{\Z}{\mbox{$\Bbb Z$}}

\newcommand{\cp}{\scriptsize \copyright}

\newtheorem{theorem}{Theorem}

\def\Aut{\mbox{Aut}}
\def\Hol{\mbox{Hol}}

\def\PGL{\mbox{PGL}}
\def\AGL{\mbox{AGL}}

\def\Gal{\mbox{Gal}}

\def\arccos{\mathop{\rm arccos}}
\def\arccosh{\mathop{\rm arccosh}}
\def\arcsinh{\mathop{\rm arcsinh}}

\def\iso{\simeq}
\def\P{\mbox{\cal P}}

\begin{document}

\maketitle

\begin{abstract}

Algebraic methods are used to construct families of
unramified abelian extensions of some families of number fields with specified
Galois groups.

\end{abstract}

\section{Abelian Extensions of Commutative Algebras}

Consider the problem of constructing split extensions 
$G\rtimes A$, with an eye to forcing the abelian kernel $A$ to
correspond to an unramified extension of an extension with group $G$,
in the case where $A$ is cyclic, so that $G \rtimes A$ is a metacyclic
group.

In the case where $G$ is cyclic of odd prime power order, we can
construct all cyclic extensions with group $G$ by the method of
\cite {Smith2}.  In the case of extensions of even degree, this
requires modification, but even so \cite {ShM} suffices to construct
these over number fields.  In any case, it is easy (\cite {Dentz}) to
see that the construction of \cite {Smith2} will give many cyclic
extensions even when it will not give all of them.

Let $n>1$ be an integer, for which we seek to construct corresponding
cyclic extensions with group $Z_n$, and $\zeta$ be primitive $n$-th
root of unity.  For $q \in \Z_{(n)}$, let $<q>$ denote the reduction
modulo $n$ to the range $0 \le <q> < n$.  Let $b_i$ denote $\phi (n)$
algebraically transcendental elements, indexed by the set of integers
$I_n$ between 0 and $n$ and prime to $n$, and let $c_i$ be such that
$c_i^n = b_i$.

Now set
$$e_j = \prod_{i,j \in I_n} c_i^{<j/i>},$$
$$r_i = \sum_{j \in I_n} e_j \zeta^{ij} , 0 \le i < n,$$
$$\P_n = \prod_{i=0}^{n-1} (x - r_i).$$

If $F$ is a field with characteristic prime to $n$ which is
disjoint from $B(\zeta)$, where $B$ is the prime field, then
by specializing $b_i$ to a conjugate orbit of values in
$F(\zeta)$, we will produce a polynomial over $F$ which
is either reducible or cyclic of degree $n$--in general, the
latter.

As explained in \cite {Smith1}, we can use this same construction to
obtain metacyclic groups which are subgroups of the holomorph
$\Hol(Z_n)$ of $Z_n$--that is to say, subgroups of the group of
invertible affine transformations 
$z \mapsto az + b$ of $Z_n$.  An automorphism of the 
cyclic group of degree $\phi(n)$ acting on the $e_i$'s will 
permute these; as a result the coefficients of the
polynomial $\P_n$ in terms of the $b_i$
are invariant under these automorphisms.

If we choose a normal basis $\omega_i$ for field $\Q(\zeta)$ over
$\Q$ (or over the prime field in nonzero characteristic), we can write
$$b_j = \sum_{i,j \in I_n} a_j \omega_{<ij>}.$$
We then have that the $a_i$ are elements of $F$, and
moreover that the coefficients of $\P_n$ when expressed
in terms of these are also invariants under automorphisms of
$Z_n$.

Now suppose that we have an extension $K/F$ of degree $\phi(n)$, 
with abelian Galois group isomorphic to the invertible elements
$(\Z/n\Z)^*$ of $\Z/n\Z$, and an explicit choice of isomorphism.  
If we have $a_i$ which are conjugate
in a manner which corresponds to the index $i$ by this isomorphism, 
then the coefficients
of the polynomial obtained by specializing $\P_n$ to these
$a_i$ will be in $F$, and so if the polynomial is irreducible, we
will obtain a an extension with Galois group $\Hol(Z_n)$.
In particular, if $n$ is an odd prime power, we have that
$(\Z/n\Z)^*$ is cyclic of degree $\phi(n)$, and we have
a metacyclic extension with group $Z_{\phi(n)} \rtimes Z_n$.

In the above case, the $a_i$ are a complete orbit of values for an
abelian extension with Galois group isomorphic to $(\Z/n\Z)^*$; and
the orbit is given by an isomorphism between $\Gal(K/F)$ and
$(\Z/n\Z)^*$.  We can generalize this to an orbit in an abelian
extension of {\em commutative algebras}.  What this means concretely
is that we have a direct sum of $l$ copies of an abelian extension of
degree $m$, such that $lm = \phi(n)$.  The Galois group of the abelian
extension of degree $m$ combines with an abelian permutation of the $l$
copies, to produce an automorphism of the algebra with the required
abelian Galois group.  The orbit of values for this algebra can now be
used for the $a_i$ just as in the case where this was a field
extension.

To give a concrete example of this, suppose that we have two copies of
the field $\Q(\sqrt{-47})$.  The Galois action 
$\sqrt{-47} \mapsto -\sqrt{-47}$ and the permutation of the two copies 
combine to give us a cyclic algebra extension of $\Q$ with Galois
group $Z_4$. A possible orbit of values in this might be
$$a_1 = -13 + \sqrt{-47}, a_2 = (-21-\sqrt{-47})/2,$$
$$a_4 =  -13-\sqrt{-47},a_3 = (-21+\sqrt{-47})/2.$$
If we now compute corresponding $b_i$ by taking
$$b_i = \sum a_j \zeta^{2^{(i+j)}},$$
then we have four roots of the polynomial 
$$x^4 - 47x^3 + 519x^2 + 47x + 1,$$
permuted in the order $x \mapsto (8x^3 -377x^2 + 4186x + 234)/5$.
If we substitute the roots of this polynomial,
{\em in this order}, for the $b_i$ in $\P_5$, we obtain
the irreducible polynomial
$$x^5 - 10x^3 -2605x^2 + 5860x + 443629.$$

This is a cyclic extension of the algebra of the $a_i$, which
means it is a cyclic extension of $\Q(\sqrt{-47})$.  It is clearly
a split extension, and we easily check (as must be the case)
that its splitting field is dihedral over $\Q$.

This is how we go about constructing split extensions with
Galois group contained in $\Aut(Z_n)$ in general.  We take
an abelian extension of the appropriate type, and use values
in this subfield to define orbits of values in the $a_i$, which
can be thought of as a single orbit in an abelian algebra.
If $n$ is the power of an odd prime, and we want to construct
{\em all} metacyclic extensions with our given extension as
subfield, we need to look at the $\phi(\phi(n))$ different
generators of $(\Z/n\Z)^*$, and allow for conjugate orbits 
produced by each of these generators.  

The example also illustrates the duality principle involved.  We began
with an orbit of values $a_i$, defined as belonging to an abelian
Galois algebra extension $K/F$ with group $\Gal(K/F) \iso \Aut(Z_n)$,
together with an explicit choice of isomorphism.  We get a
corresponding abelian Galois algebra extension $H/F$ for the values
$b_i$, which a subalgebra of $K \otimes F(\zeta)$, together with
an explicit isomorphism.  Starting from $H/F$ and the selected
isomorphism, the same diagonal subalgebra construction brings us back
to $K/F$ and our original isomorphism--these are in duality.

In the case where $n$ is an odd prime power, the duality can be
interpreted as being between a cyclic {\em field} extension
$K/F$ of degree $m$, where $m$ divides $\phi(n)$, together
with a generator of its Galois group; and another cyclic field 
extension of degree dividing $\phi(n)$ and a generator for 
{\em its} group.  In the example already given, $\Q(\sqrt{-47})$
was dual to the cyclic field of conductor 235 
previously given, together with the indicated generator for its
Galois group.  If we take the Galois group in reverse order,
so that $x \mapsto (-13x^3+612x^2-6786x-234)/5$ instead, we have
that this is now dual to $\Q(\sqrt{-235})$.  In fact,
if 
$$a_1 = (-21+\sqrt{-235})/2, a_2 = -13,$$
$$a_4 = (-21 - \sqrt{-235})/2, a_3 = -13$$
then just as before, if
$$b_i = \sum a_j \zeta^{2^{(i+j)}}$$
then we obtain roots of 
$$x^4 -47x^3 + 519x^2 + 47x + 1.$$
However, these roots now permute in the order opposite to before,
and if we substitute values in this opposite order into
$\P_5$, we obtain instead
$$x^5 - 10x^3 -2605x^2 + 5680x + 167504.$$
This is another dihedral polynomial, only now, of course,
the quadratic subfield is $\Q(\sqrt{-235})$.

We wish to find extensions of $K$ which are unramified. Since
these are extensions obtained by taking $K(\zeta)$-linear
combinations of $n$-th roots in $K(\zeta)$, we first
insure that the extensions are unramified over $K(\zeta)$.
But an $n$-th root produces an unramified extension if and only
if it is the $n$-root of something whose corresponding principal
ideal is an $n$ power as an ideal.  Since adding the $n$-th roots
of unity can cover ramification only at primes dividing $n$, 
and since if $H$ is dual to $K$ for some choice of isomorphism, it
follows that an unramified extension of $K$ obtained from
one of the polynomials $\P_n$ corresponding to this isomorphism 
must come from
$n$-th roots of values in $H$ such that
these are $n$-powers as principal ideals, and that anything
constructed in this way will be unramified outside of
primes dividing $n$.  In particular, if the $H$-values in the
orbit for the $b_i$ are units, we will obtain a corresponding extension
of $K$ unramified outside of primes dividing $n$.

Thus, for example, the two dihedral degree-five polynomials we constructed
using the quartic units given by the roots of
$$x^4 -47x^3 + 519x^2 + 47x + 1$$
over $\Q(\sqrt{-47})$ and $\Q(\sqrt{-235})$ respectively, will be
unramified outside of ramification at 5.  Hence, to see if these are
unramified, we need check only at 5.  

If we take 
$$x = (2z^4+3z^3+254z^2-3300z -9334)/64625,$$
we find that
$$x^5+x^4+x^3-x^2-2x-1$$
is a polynomial giving the same splitting field; both
polynomial and field have discriminant $47^2$.  hence this
must give an unramified cyclic extension of $\Q(\sqrt{-47})$
of degree five.

On the other hand, if we take 
$$x = (3z^4 + 44z^3 + 208z^2 -6895z -55660)/22090,$$
then
$$x^5 - 35x^3 + 50x + 20$$
has the same splitting field, and is Eisenstein and hence totally
ramified at 5.  Hence it gives us a cyclic extension of 
$\Q(\sqrt{-235})$ which is unramified outside of 5.

We can check this by comparing the ratio of the discriminant of
the degree five field with the square of the discriminant of
the quadratic subfield.  By the analysis in the previous section,
these should be equal if and only if the extension is unramified, and
we find this is the case.

\section{Families of Unramified Extensions of Cyclic Extensions}

If $q=p^m$ is an odd prime power, and if the degree $n$
of one of these polynomials divides $\phi(q)$, we may replace
$n$ of the $b_i$ in the generic polynomial for cyclic
extension of degree $q$ with conjugate values of the polynomial,
and the rest with 1.  If $k$ is of multiplicative 
order $n$ mod $q$, we may make $b_{k^i} = \sigma^i r$, where
$r$ is a root of the polynomial of degree $n$ and $\sigma$
a generator of its Galois group over $\Q(t)$.  In this
way, we obtain a polynomial with coefficients in $\Z[t]$
of degree $q$, which has a subfield of degree $n$ which is
cyclic over $\Q(t)$.  

For specializations of $t$ for which
we obtain an extension of degree $n$ (nearly all, by
Hilbert irreducibility) we then have that a root of the
specialized polynomial of degree $n$ (if the polynomial
is irreducible) will give a cyclic extension of degree $q$
of the cyclic extension of degree $n$, unramified outside
of primes dividing $n$, by the previous section.

We are now in the position to create unramified families of extensions
for degrees higher than four.  To start with, we have two degree two
polynomials $x^2-tx+1$ and $x^2-tx-1$ which give families of units in
quadratic extensions.  If $n$ is an odd prime power and if we take
$b_i$ equal to 1 in $\P_n$ except for two conjugate values for $b_1$
and $b_{n-1}$, we obtain a polynomial which gives a cyclic extension
of degree $n$ over a cyclic extension of degree $\phi(n)$, and which
is unramified except at the prime dividing $n$.

Moreover, there is an especially nice description of these extensions
in terms of Chebyshev and Lucas polynomials. Let us name the roots of
$x^2-tx+1$ as $a = ((t+\sqrt{t^2-4})/2)^{1 \over n}$ and
$b = ((t-\sqrt{t^2-4})/2)^{1 \over n}$. Then one of the roots we obtain 
by substituting $a$ for $b_1$ and $b$ for $b_{n-1}$ will be $r = r_0 = \sum_{j \in I_n} a^j b^{n-j}$.
Since $r$ is a symmetric polynomial in $a$ and $b$, homogenous of degree $n$,
$r$ can be written in terms of the two elementary symmetric polynomials, $ab =1$ and $a+b = x$. 
Hence, $r$ can be expressed as a polynomial in $x$, and so $\Q(r)$ and $\Q(x)$ give the same extension of 
$\Q$. For example, if $n=5$ we have $r = ab^4 + a^2b^3 + a^3b^2 + a^4b$, which can be written in terms of
the elementary symmetric polynomials as $ab(a+b)^3 - 2a^2b^2(a+b)$, and setting $ab=1$ and $(a+b)=x$,
this tells us that $r=x^3-2x$.

In the language of \cite {Smith3}, $x$ is the nth Chebyshev root of $t$, where we may write
$$x = ((t+\sqrt{t^2-4})/2)^{1 \over n} + ((t-\sqrt{t^2-4})/2)^{1 \over n} = 2 \cosh(\arccosh(t/2)/n) =$$
$$2 \cos(\arccos(t/2)/n) = C_{1 \over n}(t) = \sqrt[\cp n]{t}.$$ Since $x = \sqrt[\cp n]{t}$, 
$x^{\cp n} = (t^{\cp 1 \over n})^{\cp n} = t$, and the polynomial is $x^{\cp n}-t$, where $x^{\cp n}$ is the nth Chebyshev polynomial, normalized to be monic and orthogonal over [-2, 2].

In the same way, if we set $a = ((t+\sqrt{t^2+4})/2)^{1 \over n}$ and $b = ((t-\sqrt{t^2+4})/2)^{1 \over n}$, then
we may proceed as before except that now $ab = -1$. We now obtain a Lucas polynomial root;
$$x = ((t+\sqrt{t^2+4})/2)^{1 \over n} +((t-\sqrt{t^2-4})/2)^{1 \over n} = L_{1 \over n}(t),$$
and we get $L_n(x) - t$, where $L_n(x)$ is the nth Lucas polynomial.

As an algebraic function of $t$, this polynomial defines
a single real branch, which we can express by
$$L_{1 \over n}(t) = 2 \sinh (\arcsinh ({t \over 2})/2).$$
As before, however, we are really more concerned with $p$-adic
branches, and merely note that the real branch indicates we need 
not concern ourselves with ramification at the infinite place.

We now have the following:
\begin{theorem}
Upon specialization of $t$ to a value in any number-ring
$O$, both $x^{\cp n} - t$ and $L_n(x) - t$ produce an infinite number
of distinct splitting fields which are unramified over some
extension of the field of quotients of $O$.
\end{theorem}
Proof: Any value of $t$ close enough $\wp$-adically to $a^{\cp n}$
(respectively, $L_n(a)$) for all $\wp$ dividing $n$ for some $a$ in
$O$ will converge to a $\wp$-adic root for each of these $\wp$, the
only primes which can ramify.  These will create congruence conditions
for families of unramified extensions, which will give unramified
cyclic extensions of degree $n$ over an extension field of the field
of fractions of $O$ for all but a thin set of values, and hence for an
infinite set.  
$\Box$

The situation is most easily analyzed when $n$ is a power of an
odd prime, as we see in the following two theorems.

\begin{theorem}
If $p$ is an odd prime, and $q=p^k$,
let $F$ be the function field which is the 
(unique) twist of the 
$q$-th roots of unity by the algebra consisting of $\phi(q)/2$
copies of $\Q(\sqrt{t^2-4})$.  For all but a thin set of
values, specializing $t$ to a value in $j \in \Z$ produces a cyclic
extension $F_j$ of $\Q$ of degree $\phi(q)$.  If $j$ is not a $p$-th
Chebyshev power, that is, if it is not the case that 
$j = m^{\cp p}$ for some integer $m$, and if $j$ is congruent to
a fixed point $m$ of $x^{\cp p}$, (that is, a root of $x^{\cp p}-x$,) modulo
$p^{2k+1}$ if $m \equiv \pm 2 \pmod p$, or modulo $p^{k+1}$ otherwise,
then the splitting field of
$$\Q(\sqrt[\cp q]{j})$$
defines an unramified cyclic extension of $F_j$.
\end{theorem}
Proof:  We have $p-1$ $p$-adic fixed points, one for each
congruence class mod $p$, since these lift from
the factorization of $x^{\cp p}-x = x^p-x$ mod $p$.
These are also the fixed points for $x^{\cp q}$, since this is
$k$ iterations of $x^{\cp p}$.

Let $u$ be a fixed point.  Then expanding around $u$,  we have
$$x^{\cp p} = u + (x^{\cp p})'(u)(x-u) + ... .$$
Unless the fixed point is $\pm 2$, $(x^{\cp p})'(u)$ will have
$p$-adic valuation $1 \over p$
$\Box$

\begin{theorem}
If $p$ is an odd prime, $q=p^k$,
let $F$ be the function field which is the
(unique) twist of the $q$-th roots of unity by the algebra consisting
of $\phi(q)/2$ copies of $\Q(\sqrt{t^2+4})$.  For all but a thin set of
values, specializing $t$ to a value in $j \in \Z$ produces a cyclic
extension $F_j$ of $\Q$ of degree $\phi(q)$.  If $p$ is congruent to
one mod four, and if it is not the case that
$j = L_p(m)$ for some integer $m$, and if $j$ is congruent to a fixed
point of $L_p$, that is, a root of $L_p(x)-x$, modulo $p^{k+1}$, then
the splitting field of
$$\Q(L_{1 \over q}(j)).$$
defines an unramified cyclic extension of $F_j$.
Similarly, if $p$ is congruent to three mod four, and if $j$ is congruent
to $\pm \sqrt{-2}$ mod $p^{2k+1}$, or congruent to one of the other fixed
points of $L_p$ mod $p^{k+1}$, then
the splitting field of
$$\Q(L_{1 \over q}(j)).$$
defines an unramified cyclic extension of $F_j$.
\end{theorem}

We can do something similar with polynomials of higher degree giving
cyclic extensions. There are certain well-known polynomials over
$\Z[t]$ which have cyclic Galois group and norm term $\pm 1$.
We have polynomials of degrees 3, 4, and 6 such that the
indeterminate $t$ is a rational function of the roots, and
such that the roots transform by a linear rational function:
$$x^3 - tx^2 - (t+3)x - 1,$$
$$x^4 -tx^3 -6x^2 + tx +1,$$
$$x^6-2tx^5+5(t-3)x^3 + 20x^3 -5tx^2 + 2(t-3)x+1.$$

These are normalized by the condition that the polynomial
of degree $n$ is {\em apolar} to the corresponding cyclotomic
polynomial of $n$-th roots of unity, and that the indeterminate
$t$ is a positive integral multiple of the trace.  Apolarity
means that these polynomials can be written in the form
$$p(x-\omega)^n + q(x-\omega')^n,$$
where $\omega$ and $\omega'$ are the two primitive $n$-th roots
of unity.

We also have other unit-generating polynomials with cyclic
Galois groups for degrees 4, 5, and 6, whose genus is greater 
than 0:
$$x^4-t^2x^3-(t^3+2t^2+4t+2)x^2-t^2x+1,$$
$$x^5-t^2x^4-(t^3+6t^2+10t+10)x^3-(t^4+5t^3+11t^2+15t+5)x^2$$
$$+(t^3+4t^2+10t+10)x - 1,$$
$$x^6-tx^5-(t^2-5t+12)x^4+(t^3-4t^2+10t-2)x^3$$
$$-(t^3-6t^2+17t-21)x^2-(t^2-3t+6)x-1.$$
These polynomials are discussed in \cite {Emma},
\cite {Wash}, and \cite {XX}, and are related to modular functions.

For any such polynomial of degree $n$, we may produce two
corresponding polynomials of degree $\phi(q)$ so long as
$n$ divides $\phi(q)$ by the process described above.

For degree five, using the apolar unit-generating cyclic polynomial of
degree four, we obtain in one direction
$$z^5-10z^3+20z^2+(5t^2+65)z-t^3-2t^2-16t-28.$$
We obtain the same polynomial with $-t$ in the place of
$t$ from the other direction.
The Tschernhausen transformation
$$x=((t+6)z^4+(t-14)z^3+(t^2-11t-22)z^2-(3t^2+3t-150)z+4t^3+20t^2+68t-8)/D,$$
$$D=(t^3+4t^2+60t+32)$$
give us the transformed polynomial
$$x^5+10x^3-5tx^2-15x-t^2+t-16.$$
This has Galois group $Z_5$ over its cyclic quartic subfield,
given by 
$$y^4+5(t^2+16)y^2+5(t^2+16)(t-2)^2.$$

The method by which this transformation was discovered was
as follows:  for various specializations of $t$ of a polynomial of degree $n$ in $x$ over $\Z[t]$, one performs
the POLRED algorithm described in \cite {Cohen}, which gives a list of n polynomials.  One then
notes that for many of these values of $t$, we have polynomials on each list
which appear to be falling into a pattern, in the sense that they
all appear to be specializations of a polynomial over $\Z[t]$.
Interpolating these values gives us a new polynomial over
$\Z[t]$; we may then use factorization over field extensions
to both verify that the new polynomial gives the same field
extension as the old one, and to find the Tschernhausen transformation.

The implementation of POLRED on Pari and of polynomial factorization
in Maple was used to produce the above transformed polynomial as
well as the other transformed polynomials listed below.

It is evident from this that a version of POLRED which works over 
$\Z[t]$ might be possible and would be desirable.  This is interesting
also from the point of view of geometry, as the corresponding plane
algebraic curve is in part desingularized.

If $v_1$ and $v_2$ are the two roots of
$$v^2-(t^2-t+16)v+t^3+2t^2+16t+28,$$
then one root of the above degree five polynomial can be expressed as
$$L_{1 \over 5}(v_1) + L_{1 \over 5}(v_2),$$
with the other roots being nearly as easy to write in this
way. It is now possible to use this expression for the roots of in
terms of Lucas polynomial radicals to began to analyze unramified extensions.

If we substitute an integer value for $t$, if the polynomials
involved are irreducible, we will obtain an extension of a cyclic
quartic field unramified outside of 5.  We may then check locally
for where the degree five polynomial factors 5-adically, and
obtain the following:

\begin{theorem}
Let $t$ be a rational integer congruent to -5, 1, 3, or 9 mod 25.
If the polynomial 
$$y^4+5(t^2+16)y^2 + 5(t^2+16)(t-2)^2$$
is irreducible, it gives a cyclic extension of $\Q$, and if
$$p(x) = x^5+10x^3-5tx^2-15x-t^2+t-16$$
is likewise irreducible, it gives an unramified cyclic extension
of $\Q(y)$.  
\end{theorem}
Proof:
For each of these $t$ values, we can find a corresponding
$u$ such that $p(u)$ is 5-adically less than $p'(u)^2$,
and hence provides a starting value which converges to a 
root of $p(x)$ by Newton's method.  For example, if
$t=1+25s$, then $u=-9-10s$ will work.  Since we have a
5-adic root, the extension is unramified over 5, but
since it can possibly ramify only over 5, it is unramified
everywhere.
$\Box$

This polynomial over $\Q(t)$ has the advantage of being ``geometric'', so that
its splitting field has no subfields algebraic over $\Q$.  This
allows us to use it to find quadratic extensions with unramified
degree five cyclic extensions.

We obtain a $D_5$ extension when $t^2+16=5s^2$, which gives
us a Pell's equation if we wish integral values of $t$.
A recurrence relation for this is given by
$a_0=2, a_1=8, a_2=22$ and
$$a_i = 3a_{i-1}-a_{i-2}.$$
Modulo 25, this gives us values of $\pm 2, \pm 3, \pm 8$.
Since by the proof of the previous theorem, the values 
congruent to 3 will give us a 5-adic root, we obtain a
family of unramified cyclic quintic extensions of quadratic 
extensions in these cases.  If we set $b_0=-22, b_1=2728, b_2=41266478$,
and
$$b_i = 15127b_{i-1}-b_{i-2},$$
we obtain our family of unramified extensions.

We may express these values in terms of Lucas numbers.  These
are the numbers given by the recurrence relation
$L_0 = 2$, $L_1 = 1$, 
$$L_i = L_{i-1} + L_{i-2}.$$
If $\tau$ and $\bar \tau$ are the two roots of $x^2-x-1$, then
$$L_i = \tau^i + {\bar \tau}^i.$$
We assume $L_i$ is defined for all integers, positive and
negative.

In terms of Lucas numbers, we have
$$b_i = 2 L_{20i-5}.$$

Hence we obtain the following:
\begin{theorem}
Substituting a value of $t$ such that
$$t = 2 L_{20i-5}$$
into
$$x^5+10x^3-5tx^2-15x-t^2+t-16$$
leads to an unramified cyclic extension of the quadratic extension
$$\Q({\sqrt{-\sqrt{{t^2+16} \over {500}}}})$$
if the polynomial is irreducible.
Moreover, it is irreducible for
all but a finite number of integers $t$, and hence this procedure
generates an infinite family.
\end{theorem}

As an algebraic curve,
the degree five polynomial has genus one.  
Since the polynomial is of degree two in $t$, we very
nearly have a Weierstrass model already.  Solving for
$t$ in terms of $x$, and substituting
$$y^2 = 4x^3+x^2-2x-7,$$
we obtain
$$t = {{(x+3)y-5x^2+1} \over 2}$$
in terms of $x$ and $y$ on the Weierstrass model
$$y^2 = 4x^3+x^2-2x-7.$$

By Mordell's theorem, this has only a finite number of integral
points. Moreover, the curve (of conductor 50) has a reduced minimal
model $$y^2+xy+y = x^3-x-2,$$ which by the tables of \cite {Cre} has
rank 0 and a torsion subgroup with three elements, the nonzero
elements of which correspond to $t=-22, x = 2$, and $t=3, x=2$.  Hence
in this case we can completely analyze the function field from the
point of view of irreducibility under specialization, and conclude that
it is irreducible for all rational specializations excepting $t=-22, 3$.

From the genus one unit-generating cyclic quartic of \cite {Wash},
we obtain in one direction the polynomial
$$z^5-10z^3+5(t^3+2t^2+4t+4)z^2-5(t^4 +2t^3+4t^2+8t+3)z$$
$$+t^7+4t^6+10t^5+22t^4+29t^3+26t^2+20t+4.$$
We may transform this as before, obtaining
$$x^5-10x^3-5x^2t^2+5(t^3+2t^2+4t+5)x-(t^3+2t^2+5t+8)t,$$
and in any case, we have again a polynomial with Galois
group $F_{20}$, and in this case, the cyclic quartic 
subfield is given by
$$y^4+5(t+2)(t^2+4)ty^2+5(t^2+4)(t+2)^2(t-1)^2t^2.$$
Because of the modular interpretation of the unit-generating
polynomial, we have a modular interpretation of this 
$F_{20}$ polynomial, which because of its simple form, may be 
a rather natural one.

As before we may prove this:

\begin{theorem}
Let $t$ be a rational integer congruent to -1 mod 5, 
-8 or -2 mod 25, or 0 mod 125.
If the polynomial 
$$y^4+5(t+2)(t^2+4)ty^2+5(t^2+4)(t+2)^2(t-1)^2t^2.$$
is irreducible, it gives a cyclic extension of $\Q$, and if
$$(x^2-5)^2x+(20x-8)t-5(x-1)^2t^2+(5x-2)t^3-t^4,$$
is likewise irreducible, it gives an unramified cyclic extension
of $\Q(y)$.
\end{theorem}
Proof: As before.
$\Box$

We have a quadratic rather than a quartic subfield when
$t^2+4=5s^2$, and this Pell's equation has solutions
$a_0=1, a_1=4, a_2=11$ and $a_i=3a_{i-1}-a_{i-2}$.  Modulo
5, this alternates between 1 and $-1$, and we obtain
the value congruent to $-1$ mod 5 by
$b_0=-1, b_1=4, b_2=29$ and $b_i = 7b_{i-1}-b_{i-2}$. This
we may express in terms of Lucas numbers by
$$b_i = L_{4i-1}.$$
As before, we have thus a family of cyclic quintic 
extensions of quadratic extensions; and also
as before, this is easily seen to be an infinite family; geometrically,
we have a curve of genus three, and hence only finitely
many rational points on it.

\begin{theorem}
Substituting a value of $t$ such that
$$t = L_{4i-1}$$
into
$$x^5-10x^3-5x^2t^2+5(t^3+2t^2+4t+5)x-(t^3+2t^2+5t+8)t$$
$$(x^2-5)^2x+(20x-8)t-5(x-1)^2t^2+(5x-2)t^3-t^4$$
leads to an unramified cyclic extension of the quadratic extension
$$\Q(\sqrt{-t(t+2)\sqrt{5t^2+20}})$$
if the polynomial is irreducible.
Moreover, it is irreducible for
all but a finite number of integers $t$, and hence this procedure
generates an infinite family.
\end{theorem}

\begin{theorem}
Let $t$ be a rational integer congruent to 7 or 11 mod 25,
-2 or 0 mod 125, or 989 mod 3125.
If the polynomial 
$$y^4+5(t+2)(t^2+4)ty^2+5(t^2+4)(t+2)^2(t+1)^2t^2$$
is irreducible, it gives a cyclic extension of $\Q$, and if
$$(x+4)(x-1)^4+20(x-1)^2t+(10x^2-20x+26)t^2+(5x^2-10x+13)t^3+(-5x+6)t^4+2t^5$$
is likewise irreducible, it gives an unramified cyclic extension
of $\Q(y)$.
\end{theorem}
Proof: As before.
$\Box$

\begin{theorem}
Substituting a value of $t$ such that
$$t = L_{20i-15}$$
or
$$t = L_{100i-25}$$
into
$$(x+4)(x-1)^4+20(x-1)^2t+(10x^2-20x+26)t^2+(5x^2-10x+13)t^3+(-5x+6)t^4+2t^5$$
leads to an unramified cyclic extension of the quadratic extension
$$\Q({\sqrt{-t(t+2)\sqrt{{t^2+4} \over {125}}}}).$$
if the polynomial is irreducible.
Moreover, these are irreducible for
all but a finite number of integers $t$, and hence each of these procedures
generates an infinite family.
\end{theorem}

Let us turn now to consideration of polynomials of degree seven 
with splitting fields which have Galois group
$F_{42}$, and which are unramified over their cyclic subfields
of degree six.  In considering how to construct these from unit-generating
polynomials, it is useful to consider first the subfields of
the unit-generating polynomials.  The apolar unit-generating
polynomial has a cubic subfield given by 
$$x^3+tx^2+(t-3)x-1$$
and quadratic subfield
$$\Q(\sqrt{t^2-3t+9}).$$
The other unit-generating polynomial, of genus two, has the
same cubic subfield, but has a quadratic subfield given by
$$\Q(\sqrt{(t-2)^2+4}),$$
that is, by the roots of $x^2-(t-2)x-1$.

If we use the apolar unit-generating polynomial of degree six to
construct an $F_{42}$ polynomial of degree seven, we obtain
in one direction
$$x^7-21x^5+70x^4-105x^3-28(4t^2-12t+33)x^2+7(96t^2-288t+859)x$$
$$+64t^3-1264t^2+3792t-9642.$$
The polynomial we obtain in the other direction can be obtained
from this one by the substitution $t \mapsto 3-t$.

The cubic and quadratic subfields of this can be
obtained by twisting the cubic and quadratic subfields of
the unit-generating polynomial by the roots of a polynomial
giving the cyclic cubic field of conductor 7 and by
$\sqrt{-7}$ respectively.  By  a twist of a cyclic extension of
degree $n$ by another of degree $n$, I mean a cyclic extension
of degree $n$ contained in the compositum, which is not identical
to either of the original extensions.

This gives us the following theorem:

\begin{theorem}
Let $t$ be a rational integer congruent to 0, 5, 8, 17, 20, or
25 mod 49.
If the polynomial 
$$y^3-ty^2-(2t^2-7t+21)x+t^3+28$$
is irreducible, it gives a cyclic extension of $\Q$, and if
$$x^7-21x^5+70x^4-105x^3-28(4t^2-12t+33)x^2+7(96t^2-288t+859)x$$
$$+64t^3-1264t^2+3792t-9642.$$
is likewise irreducible, it gives an unramified cyclic extension
of $\Q(\sqrt{-7(t^2-3t+9)}, y)$.
\end{theorem}
Proof: As before.
$\Box$

If we use the genus-two unit-generating polynomial, we obtain
in one direction a polynomial which we may simplify by transformation
and so obtain the following theorem:

\begin{theorem}
Let $t$ be a rational integer congruent to 2 mod 7, or to
$-21,-18,-16,-8,11$ mod 49, or 743 mod 2401.
If the polynomial 
$$y^3-ty^2-(2t^2-7t+21)x+t^3+28$$
is irreducible, it gives a cyclic extension of $\Q$, and if
$$x^7+21x^5-7(t^2-4t+10)x^4+28(t^2-3t+15)x^3-7(5t^3-8t^2+12t+72)x^2+$$
$$7(5t+6)(2t^2-7t+22)x-t^5-20t^4-94t^3+410t^2-1584t+1224$$
is likewise irreducible, it gives an unramified cyclic extension
of $\Q(\sqrt{-7(t-2)^2-28}, y)$.
\end{theorem}
Proof: As before.
$\Box$

Going in the other direction, we will obtain a $F_{42}$ extension
with the same quadratic subfield.  However, it seems to me more
interesting to make the substitution $t \mapsto 3-t$ after doing
this, and obtain another $F_{42}$ extension with the same cyclic
cubic subfield.  Doing this and making another simplifying
transformation gives us:

\begin{theorem}
Let $t$ be a rational integer congruent to 3 mod 7, or to
$-17,-5,5,7,13$ mod 49, or 743 mod 2401.
If the polynomial 
$$y^3-ty^2-(2t^2-7t+21)x+t^3+28$$
is irreducible, it gives a cyclic extension of $\Q$, and if
$$x^7-7tx^5-7(t^2-4t+11)x^4+28(t^2-t+3)x^3+7(3t^2-13t+36)tx^2+$$
$$7(t^4-18t^3+68t^2-176t+192)x-t^5-23t^4+184t^3-816t^2+1536t-2304$$
is likewise irreducible, it gives an unramified cyclic extension
of $\Q(\sqrt{-7(t-1)^2-28}, y)$.
\end{theorem}
Proof: As before.
$\Box$

We may also use the unit-generating polynomials of degree six to
produce polynomials of degree nine with Galois group $F_{54}$.
When we do this, we discover that using the genus zero polynomial
leads to a reducible polynomial, and that using the genus two
polynomial leads to two different extensions.  These were not reduced
by the POLRED method, which failed to produce the requisite family of polynomials
in $\Z[x]$.

\begin{theorem}
Let $t$ be a rational integer congruent to $-12,11,6$ or 11 mod 27.
If the polynomial 
$$y^3-3(t^2-3t+9)x+(t-6)(t^2-3t+9)$$
is irreducible, it gives a cyclic extension of $\Q$, and if
$$x^9+27x^7-9(t^3-4t^2+11t-6)x^6+27(t-2)(t^2-4t+8)(t^2-3t+9)x^5-$$
$$9(t^2-4t+8)(t^2-3t+9)(t^4-7t^3+26t^2-48t+36)x^4-$$
$$3(7t^2-18t+63)(t^2-3t+9)(t^2-4t+8)^2x^3-$$
$$27(t^2-5t+10)(t^2-4t+8)^2(t^2-3t+9)^2x^2+$$
$$9(t^3-6t^2+18t-24)(t^2-3t+9)^2(t^2-4t+8)^3x-$$
$$(t^6-11t^5+61t^4-213t^3+475t^2-660t+468)(t^2-3t+9)^2(t^2-4t+8)^3$$
is likewise irreducible, it gives an unramified cyclic extension
of $\Q(\sqrt{-3(t-2)^2-12}, y)$.
\end{theorem}
Proof: As before.
$\Box$

Once again, under the assumption that getting the cubic subfields to
agree is marginally more interesting than getting the quadratic
subfields to agree, I make the substitution of $3-t$ for $t$ when
going in the other direction:

\begin{theorem}
Let $t$ be a rational integer congruent to 5 or 10 mod 27.
If the polynomial 
$$y^3-3(t^2-3t+9)x+(t-6)(t^2-3t+9)$$
is irreducible, it gives a cyclic extension of $\Q$, and if
$$x^9+27x^7+9(t^3-5t^2+14t-18)x^6-27(t-1)(t^2-2*t+5)(t^2-3t+9)x^5+$$
$$9(t^2-2t+5)(t^2-3t+9)(t^3-2t^2+5t+12)x^4-$$
$$3(t^2-3t+9)(10t^2-33t+90)(t^2-2t+5)^2x^3+$$
$$27(t^2-3t+6)(t^2-2t+5)^2(t^2-3t+9)^2x^2-$$
$$9(t^2-2t+5)^3(t^2-3t+9)^3x+$$
$$(t^4-5t^3+27t^2-54t+135)(t^2-3t+9)^2(t^2-2t+5)^3$$
is likewise irreducible, it gives an unramified cyclic extension
of $\Q(\sqrt{-3(t-1)^2-12}, y)$.
\end{theorem}
Proof: As before.
$\Box$

This ends the list of geometric extensions producing families of
unramified cyclic extensions of cyclic extensions which we can
produce by this method.  There are, however, an infinity of
non-geometric possibilities as well; for example,
consider the extension of degree seven with Galois group
$F_{42}$ obtained from the unit-generating polynomial of
degree three $x^3+tx^2+(t-3)x-1$.  Transforming this via
the POLRED method and proceeding as before, we obtain:

\begin{theorem}
Let $t$ be a rational integer congruent to $-16,-11,-5,0,6$
or 11 mod 49, or 743 mod 2401.
If the polynomial 
$$y^3-ty^2-(2t^2-7t+21)x+t^3+28$$
is irreducible, it gives a cyclic extension of $\Q$, and if
$$x^7-14x^4-7(t-3)x^3+14tx^2-28x+t^2-11t+33$$
is likewise irreducible, it gives an unramified cyclic extension
of $\Q(\sqrt{-7}, y)$.
\end{theorem}
Proof: As before.
$\Box$

\section{Unramified 2-elementary extensions of $\PGL_3(2)$}

In \cite {Smith1}, I explain a method for constructing split
extensions with Galois groups a subgroup of the holomorphs of
2-elementary abelian groups.  In the case of the holomorph
$\Hol (2^3)$, this means a split extension of $\PGL_3(2)$, called
among other things $\AGL_3(2)$, $2^3.\PGL_3(2)$, or 8T48.

Back in the last millenium the author constructed a polynomial of degree eight over
$\Q(b_1, \cdots, b_7)$ using seven indeterminates $b_i$,
which had coefficients which were resolvents for 
the roots of a polynomial with Galois group 
$\PGL_3(2)$. When the seven roots of a polynomial of degree seven were properly ordered 
and substituted into the degree eight polynomial, one then obtained a split extension of
$\PGL_3(2)$ so long as the roots were not squares of elements 
in the splitting field, and generated the splitting field.

One especially interesting polynomial for this purpose is
$$x^7+(3u-2)x^6+(2u^2+2u-3)x^5+(3u^2-3u+t)x^4-(u^3-u^2+2u-t)x^3$$
$$-u^2(u+4)x^2+u^2(u-3)x+u^3.$$
This is a transformation the author obtained of the LaMacchia polynomial \cite {LaM}, 
which has Galois group $\PGL_3(2)$ over $\Q(t,u)$.  It has the
particularly interesting (for us) property that the norm
term is simply a power of one of the coefficients.  This meant
that when the roots were property ordered and substituted for the
$b_i$ in the mystery degree eight polynomial, I was able to obtain a polynomial with
Galois group $\Hol (2^3)$ which is very well suited for finding
unramified extensions; this polynomial is
$$x^8+(-20u^2-4t+8u)x^6-16(u-3)u^2x^5$$
$$+(16u^5-42u^4+36u^2t+80u^3+6t^2-24ut-48u^2)x^4$$
$$-32(7u^3+ut-7u^2+3t-6u)u^2x^3$$
$$+(-96u^7-32tu^5+588u^6-60tu^4+1288u^5-12t^2u^2-128tu^3+$$
$$64u^4-4t^3+24t^2u+96u^2t-256u^3)x^2$$
$$+16(24u^6+79u^5-2tu^3+59u^4+3t^2u+14u^2t-16u^3+3t^2-12ut-96u^2)u^2x$$
$$-48t^2u^2+1024u^4-1024u^5-164tu^6+256tu^3+54t^2u^4-8t^3u+48t^2u^3$$
$$-2624u^6-352u^7+t^4+753u^8+400u^9+64u^10-4t^3u^2+16t^2u^5-392tu^5-32tu^7.$$

As one particularly interesting special case, if we set $u=-1, t= 4t+2$ in the
above polynomials, and transform the second, we obtain
$$x^7-5x^6-3x^5+8x^4+6x^3-3x^2-4x-1+4tx^3(x+1)$$
as a polynomial with Galois group $\PGL_3(2)$ over
$\Q(t)$, and
$$z^8-(36+16t)z^6+64z^5+(96t^2+336t-42)z^4+(128-256t)z^3$$
$$-(256t^3+960t^2-16t+68)z^2+(1792t-320)z+(256t^4+768t^3-160t^2-592t+17)$$
with Galois group $\Hol (2^3)$ over $\Q (t)$.  
Transforming  this by the POLRED method, we obtain
$$(y+1)(y^7-y^6-11y^5+y^4+41y^3+25y^2-34y-29)-t(2y+3)^2.$$

The polynomial discriminant of the degree seven polynomial in $x$ and the degree 
eight polynomial in $y$ is exactly the same,
$$(6912t^4-3456t^3-95472t^2+23976t-1417)^2.$$
This leaves little, if any, possibility that the splitting field
for the degree seven polynomial is not a subfield of the splitting
field for the degree eight polynomial; however since the mystery polynomial
originally used to construct these has been lost, I will point out 
that any vestige of doubt is removed by the fact that the a root of the degree seven
polynomial can be expressed in terms of a root of the degree 28 polynomial which
is the polynomial satisfied by the sum of two distict roots of the degree eight
polynomial. Maple can do the necessary computations without difficulty, but
the results are so unwieldy I do not give them here.

From the manner of its construction, which involved adjoining square roots of units
in the splitting field for the degree seven polynomial, the only possible ramification
is at 2 or infinity, and in any case for any integer value of $t$, 
$6912t^4-3456t^3-95472t^2+23976t-1417$ is odd. The stem field defined by specializing $t$
to an integer in the degree seven polynomial has three real and four complex embeddings,
and for the degree eight polynomial we get four real and four complex embeddings, so we have no
ramification at any infinite place. We can also directly check the ramification at 2; the degree seven polynomial is 
identically $x^7+x^6+x^5+x^2+1$ for any integer specialization of $t$, and the degree eight
polynomial factors as $y(y^7+y+1)$ for odd $t$ and $(y+1)(y^7+y^6+y^5+y^4+y^3+y^2+1)$ for even
$t$, so there is no ramification over 2. Odd primes $p$ ramify (tamely) in the degree seven stem 
field if and only if $p$ divides $6912t^4-3456t^3-95472t^2+23976t-1417$, and similarly for 
the degree eight stem field if the polynomial is irreducible.

We therefore have this:
\begin{theorem}
Let $t$ be a rational integer.
For all but a finite number of exceptional $t$, the roots of
$$x^7-5x^6-3x^5+8x^4+6x^3-3x^2-4x-1+4tx^3(x+1)$$
gives a $\PGL_3(2)$ extension of $\Q$, and the roots of
$$(y+1)(y^7-y^6-11y^5+y^4+41y^3+25y^2-34y-29)-t(2y+3)^2.$$
gives an unramified $2^3$-elementary extension of the $\PGL_3(2)$ extension.
\end{theorem}
Proof: By Hilbert irreducibility and the above.
$\Box$

As an immediate result, we have constructed an infinity of
$\PGL_3(2)$ Galois extensions of $\Q$ such that the Galois
group acts faithfully on a quotient of the class group, in just
the way that three-by-three square matrices over the field of two
elements act on three-vectors in such a field.

In \cite {Malle}, Gunter Malle gives two polynomials with Galois group 
$\PGL_3(2)$ over three indeterminates and one over four indeterminates. 
It seems likely that many more families of unramified $2^3.\PGL_3(2)$
extensions could be constructed using these. Malle also gives a two-parameter
family of $2^3.\PGL_3(2)$ extensions, and mentions the problem of finding
a $2^3.\PGL_3(2)$ polynomial whose discriminant is the square of an odd prime,
and provides an example. Since $6912t^4-3456t^3-95472t^2+23976t-1417$ seems to
be fecund as a prime-generating polynomial, many more such polynomials
are easy to provide; by Bunyakovsky's conjecture we of course expect the
number of such primes to be infinite.

\end{document}